\RequirePackage{amsmath}
\documentclass[a4paper,12pt]{article}%

\usepackage{mathptmx}
\usepackage[english]{babel}

\usepackage{amsfonts}
\usepackage{amssymb}
\usepackage{amsthm}
\usepackage{graphicx}
\usepackage{pdfpages}
\usepackage{relsize}

\newcommand{\setx}[1]{ \{ #1 \} }
\newcommand{\setminx}[1]{ \setminus \{ #1 \} }

\newcommand{\mytext}[1]{ \: \textrm{#1} \: }
\newcommand{\mysetdescr}[2]{\left\{ #1 \: \left| \: #2 \right. \right\} }
\newcommand{\darr}{{\downarrow \,}}
\newcommand{\uarr}{{\uparrow \,}}
\newcommand{\ouarr}{{\uparrow_{_{_{\!\!\!\!\!\!\circ}}}}}\newcommand{\odarr}{{\downarrow^{^{\!\!\!\!\!\!\circ}}}}

\newcommand{\myN}{\mathbb{N}}

\newcommand{\myNkz}[1]{\setx{ 0, \ldots , #1 }}
\newcommand{\nz}{\setx{0,1,2}}

\def\rarr{\rightarrow}

\newcommand{\mf}[1]{\mathfrak{ #1 }}

\def\fC{\mf{C}}

\def\fN{\mf{N}}
\def\fS{\mf{S}}

\def\fNS{\mf{N}^{Seg}}
\def\fNL{\mf{N}^{LS}}
\def\fNU{\mf{N}^{US}}

\def\BP{\begin{proof}}
\def\EP{\end{proof}}

\newcommand{\celaga}[2]{c_{\lambda(#1),\gamma(#2)}}

\begin{document}

\theoremstyle{thmstyleone}
\newtheorem{theoremaC}{Theorem}[section]
\newtheorem{definitionaC}[theoremaC]{Definition}
\newtheorem{corollaryaC}[theoremaC]{Corollary}
\newtheorem{lemmaaC}[theoremaC]{Lemma}
\newtheorem{propositionaC}[theoremaC]{Proposition}
\newtheorem{criteriaaC}[theoremaC]{Criteria}

\renewcommand{\thetheoremaC}{\arabic{section}.\arabic{theoremaC}}





\title{A recursive approach for the determination of the nice sections of width three which have a 4-crown stack as retract \\ {\normalsize Version 9-3}}
\author{Frank a Campo}

\maketitle
\begin{abstract}
The characterization of the finite minimal automorphic posets of width three is still an open problem. Niederle has shown that this task can be reduced to the characterization of the nice sections of width three which have a non-trivial tower of nice sections as retract. In our article, we develop a recursive approach for the determination of those nice sections of width three which have a 4-crown stack as retract. We apply the approach on a sub-class $\fN_2$ of nice sections of width three and determine all posets in $\fN_2$ with height up to six which have a 4-crown stack as retract. For each integer $n \geq 2$, the class $\fN_2$ contains $2^{n-2}$ different isomorphism types of posets of height $n$.
\newline

\end{abstract}

\section{Introduction} \label{sec_Intro}

We call a finite poset $P$ {\em automorphic} if it has a fixed point free automorphism, and we call it {\em minimal automorphic} if additionally every proper retract of it has the fixed point property. Minimal automorphic posets have been present in the literature since the late eighties; an overview is provided by \cite{Schroeder_2022_MASoC}. In particular, they have found interest as forbidden retracts because a poset does not have the fixed point property iff it has a minimal automorphic poset as retract \cite[Th.\ 4.8]{Schroeder_2016}.

However, the description or characterization of minimal automorphic posets is still an open field. (Unless otherwise stated, ``poset'' always means ``finite poset'' in this article.) Since the pioneering work of Rival \cite{Rival_1976,Rival_1982}, we know that a poset of height one does not have the fixed point property iff it contains a crown, because a crown of minimal length will be a retract of it. Brualdi and Dias Da Silva \cite{Brualdi_DdaSilva_1997} generalized this result by showing that a poset does not have the fixed point property iff it has a ``generalized crown'' as retract. With respect to posets of small width, Fofanova and Rutkowski \cite{Fofanova_Rutkowski_1987} showed that a poset of width two does not have the fixed point property iff it has a 4-crown stack as retract. (The result extends even to chain-complete infinite posets.) Niederle \cite{Niederle_1989} proved in 1989 that a poset of width three does not have the fixed point property iff it has a ``tower of nice sections'' as non-trivial retract. Later on, in 2008, he extended his result to posets of width four \cite{Niederle_2008}.

Unfortunately, it is quite difficult to break down Niederle's abstract characterization to concrete posets or classes of posets. As far as we know, the only real success in this field is a result of Farley \cite{Farley_1997} from 1997, who showed that the ranked nice sections of width three are exactly the 6-crown stacks, and that a 6-crown stack is minimal automorphic iff its height $h$ is not a multiple of three. Otherwise, it has a 4-crown stack of height $\frac{2}{3} h$ as retract. Schr\"{o}der \cite[Lemma 2.6]{Schroeder_2022_MASoC} generalized one direction of this result by showing that a $2n$-crown stack of height $nm$ has a 4-crown stack of height $2m$ as retract. Moreover, in the case of $n \geq 3$, a $2n$-crown stack of height $h$ is not minimal automorphic if $h$ and $n$ have a non-trivial common factor \cite[Prop.\ 2.7]{Schroeder_2022_MASoC}.

Due to the results of Niederle \cite{Niederle_1989}, the study of the minimal automorphic posets of width three can be restricted to the investigation of the nice sections of width three. In our article, we develop and apply a recursive approach for the determination of those nice sections of width three which have a 4-crown stack as retract. After the preparatory Section \ref{sec_preparation}, we deal with the classes $\fN$ of nice sections and $\fC$ of {\em crowned sections} in Section \ref{sec_Sections}; $\fN$ is a proper subclass of $\fC$. We see that a poset $P \in \fC$ has a 4-crown stack as retract iff it has such a retract which contains a {\em split level}, and we prove several results about the existence of such levels. We see in Lemma \ref{lemma_matching} that a split level of $R$ corresponds in $P$ to a certain decomposition by a {\em retractive split}, and based on this observation we develop at the end of Section \ref{sec_Sections} a recursive approach for the determination of the nice sections with a 4-crown stack as retract.

In order to test the approach, we introduce the class $\fN_2$ of nice sections with {\em horizon two} in Section \ref{sec_horizonTwo}. For each integer $n \geq 2$, this class contains $2^{n-2}$ isomorphism types of posets of height $n$; the isomorphism types are described in Theorem \ref{theo_isomorphism} which generalizes a result of Farley \cite[Prop.\ 4.1]{Farley_1997}. We characterize in Theorem \ref{theo_RST} the retractive splits of posets $P \in \fN_2$ which yield a 4-crown stack as retract. Using this result, we determine in Section \ref{sec_Application} by means of the recursive approach all posets $P \in \fN_2$ of height up to six with a 4-crown stack as retract.

The characterization of the minimal automorphic posets of width three is regarded as a difficult task \cite[p.\ 98]{Schroeder_2016}, in parts even as ``intractable'' \cite[p.\ 126]{Farley_1997}. Also the sparse sequence of publications indicates difficulties. Nevertheless, it gives some confidence that the nice section of width three used by Schr\"{o}der \cite[p.\ 99]{Schroeder_2016} for the illustration of the problems and considered as ``nasty'' by Farley \cite[p.\ 136]{Farley_1997} is handled without a problem by the machinery developed in this article. It is the poset coded as 1011 in Section \ref{sec_Application} and shown in Figure \ref{fig_PosetsHeight4}.

\section{Preparation} \label{sec_preparation}

In this section, we introduce our notation and recapitulate definitions of structures which are in the focus of our investigation. For all other terms of order theory, the reader is referred to standard textbooks as \cite{Schroeder_2016}. For integers $i \leq j$, we write $[i,j]$ for the interval $\setx{i, \ldots , j}$.

Let $P = (X, \leq)$ be a finite poset. (We deal with finite posets only in this paper.) For $Y \subseteq X$, the {\em induced sub-poset} $P \vert_Y$ of $P$ is $\left( Y, {\leq} \cap (Y \times Y) \right)$. To simplify notation, we identify a subset $Y \subseteq X$ with the poset $P \vert_Y$ induced by it.

For $y \in X$, we define the {\em down-set} and {\em up-set induced by $y$} as
\begin{align*}
\darr y & := \mysetdescr{ x \in X }{ x \leq y }, \quad \odarr y := ( \darr y ) \setminus \setx{y}, 
\\
\uarr y & := \mysetdescr{ x \in X }{ y \leq x }, \quad \ouarr y := ( \uarr y ) \setminus \setx{y}. 
\end{align*}
Given two points $x < y$ in $P$, the point $x$ is called a {\em lower cover} of $y$ and $y$ an {\em upper cover} of $x$ if $x \leq z \leq y$ implies $z \in \setx{x,y}$ for all $z \in P$; this relation is denoted by $x \lessdot y$. A point in $P$ is called {\em irreducible} if it has a single lower cover or a single upper cover.

Let $A, B \subseteq P$. We write $A < B$ if $a < b$ for all $a \in A$, $b \in B$. If $A = \setx{a}$ is a singleton, we simply write $a < B$ and $B < a$. The notation for the relation $\leq$ is analogous.

The {\em length} of a chain is its cardinality minus 1, and the {\em height} of a poset is the maximal length of a chain contained in it. For a poset $P$, we denote its height by $h_P$. The {\em width} of a poset is the largest size of an antichain contained in it.

For a poset $P$, the {\em level sets} $P(k)$, $k \in [0,h_P]$, are recursively defined by
\begin{align*}
P(0) & := \min P, \\
P(k+1) & := \min \left( P \setminus \cup_{i=0}^k P(i) \right) \quad \mytext{for all } 0 \leq k < h_P.
\end{align*}
Additionally, we write for all $k, \ell \in [0,h_P]$ with $k \leq \ell$,
$$
P(k,\ell) := P(k) \cup P(\ell) \quad \mytext{and} \quad
P(k \rightarrow \ell ) := \cup_{i=k}^\ell P(i).
$$
For level indices $k_1 \leq k_2$ of $P$, we call the sub-poset $P(k_1 \rarr k_2)$ a {\em segment} of $P$, a {\em lower segment} if $k_1 = 0$, and an {\em upper segment} if $k_2 = h_P$.

For $n \in \myN \setminus \setx{1}$, a {\em $2n$-crown} is a poset with $2n$ points $x_{i,j}$, $i \in [0,n-1]$, $j \in \setx{0,1}$, in which $x_{0,0} < x_{0,1} > x_{1,0} < \cdots > x_{n-1,0} < x_{n-1,1} > x_{0,0}$ are the only comparabilities between different points. We call a poset $P$ of height at least 1 a {\em $2n$-crown stack} if $P(k, k+1)$ is a $2n$-crown for all $k \in [0, h_P-1]$ and all other comparabilites are induced by transitivity.

\begin{figure}
\begin{center}
\includegraphics[trim = 70 720 205 75, clip]{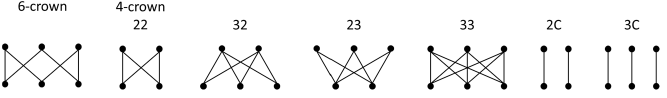}
\caption{\label{fig_Types_nm} The level pairs $P(k,\ell)$, $k < \ell$, which are possible in an automorphic poset $P$ of width at most three.}
\end{center}
\end{figure}

Given disjoint posets $P = (X, \leq_P)$ and $Q = (Y, \leq_Q)$, their {\em ordinal sum} is the poset $P \oplus Q := (X \cup Y, \leq_P \cup \leq_Q \cup \; (X \times Y) )$. For $P$ and $Q$ being disjoint antichains of size 2 or 3, we simply write ``$\#P \#Q$'' as shortcut for the isomorphism type of $P \oplus Q$, e.g., $22$ indicates a 4-crown. Additionally, we write $2C$ and $3C$ for a poset which consists of two (three) disjoint chains of height 1 without any additional comparabilities between different points. The shortcuts are illustrated in Figure \ref{fig_Types_nm}.

For a mapping $f : X \rarr Y$, we denote by $f \vert_Z$ the pre-restriction of $f$ to a subset $Z$ of its domain and by $f \vert^Z$ its post-restriction to a subset $Z$ of its codomain with $f[X] \subseteq Z$. For $X' \supseteq X$, we call a mapping $f' : X' \rarr Y$ an {\em extension} of $f$, if $f' \vert_X = f$.

An order homomorphism $r : P \rightarrow P$ is called a {\em retraction} of the poset $P$ if $r$ is idempotent, and an induced sub-poset $R$ of $P$ is called a {\em retract of $P$} if a retraction $r : P \rightarrow P$ exists with $r[X] = R$. For the sake of simplicity, we always identify $r$ with its post-restriction and write $r : P \rightarrow R$. A poset $P$ has the fixed point property iff every retract of $P$ has the fixed point property \cite[Th.\ 4.8]{Schroeder_2016}. We need the following well-known result about the existence of {\em spanning retracts}:

\begin{lemmaaC}[{\cite[p.\ 232]{Duffus_etal_1980_DPR}}] \label{lemma_R0P0}
Let $P$ be be a finite poset and $r : P \rarr R$ a retraction. There exists a retraction $v : P \rarr V$ with $V \simeq R$, $\min V \subseteq \min P$, and $\max V \subseteq \max P$.
\end{lemmaaC}

A poset is called {\em automorphic} if it has a fixed point free automorphism, and it is called {\em minimal automorphic} if additionally every proper retract of it has the fixed point property. According to \cite[Prop.\ 4.38]{Schroeder_2016}, a poset $P$ of width at most three is automorphic iff for all level indices $0 \leq k < \ell \leq h_P$, the level-pair $P(k,\ell)$ is one of the seven posets shown in Figure \ref{fig_Types_nm}. The result goes back to \cite[Th.\ 10]{Niederle_2008}.

\section{Nice sections and crowned sections} \label{sec_Sections}

\subsection{Nice sections} \label{subsec_niceSections}

The following definition is due to Niederle \cite{Niederle_1989}:

\begin{definitionaC} \label{def_Section}
A {\em section} is a two-element antichain or a poset $P$ of heigth at least one with carrier $X := \mysetdescr{ c_{k,j} } { k \in [0,h_P], j \in \nz }$ for which
\begin{itemize}
\item $C_j := c_{0,j} < \ldots < c_{h_P,j}$ is a chain for all $j \in \nz$;
\item $\setx{ c_{k,0}, c_{k,1}, c_{k,2} }$ is an antichain for all $k \in [0,h_P]$;
\item The mapping $c_{k,j} \mapsto c_{k,(j+1) \! \! \! \mod 3} \;$ for all $k \in [0, h_P]$, $j \in \nz$, is an automorphism of $P$;
\item $P(k,k+1) \not\simeq 33$ for all level indices $k \in [0, h_P-1]$.
\end{itemize}
A section $P$ is called {\em nice} if $P$ does not contain an irreducible point. The symbol $\fN$ denotes the class of nice sections. A {\em tower of sections} is an ordinal sum of sections.
\end{definitionaC}

Simple examples of nice sections are the 6-crown stacks, whereas a 4-crown stack is a tower of nice sections, i.e., of 2-antichains. 

For a section $P$ of width 3, the first and second property imply $c_{0,j} \lessdot \ldots \lessdot c_{h_P,j}$ for all $j \in \nz$. We call the chains $C_0, C_1, C_2$ the {\em main chains} of $P$. Clearly, $P(k) = \setx{ c_{k,0}, c_{k,1}, c_{k,2} }$ for all $k \in [0,h_P]$. We define mappings $\lambda : P \rightarrow [0,h_P]$ and $\gamma : P \rarr \setx{0,1,2}$ by setting $\lambda(c_{k,j}) := k$ and $\gamma(c_{k,j}) := j$, hence $\setx{p} = P(\lambda(p)) \cap C_{\gamma(p)}$ for all $p \in P$. Two consecutive levels of $P$ form a 6-crown or are of type $3C$. The level-index function $\lambda$ and the chain-index function $\gamma$ always refer to the levels and main chains of $P$; for a retract $R$ of $P$, we thus have $\setx{x} = P(\lambda(x)) \cap C_{\gamma(x)}$ for all $x \in R$.

Niederle \cite[p.\ 121 and p.\ 125]{Niederle_1989} has shown that the decomposition of a tower of sections into sections is unique, and that a poset of width at most three has not the fixed point property iff it has a proper retract which is a tower of nice sections. We conclude that a minimal automorphic poset $P$ of width at most three must be a tower of nice sections. Moreover, if $P = S_1 \oplus \cdots \oplus S_n$ is a tower of nice sections $S_1, \ldots S_n$, it is not hard to see that $P$ is minimal automorphic iff $S_i$ is minimal automorphic for every $i \in [1,n]$. The objects we have to deal with in investigating minimal automorphic posets of width three are thus the nice sections of width three.

\begin{figure}
\begin{center}
\includegraphics[trim = 70 690 370 70, clip]{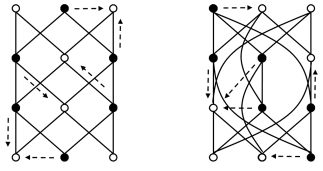}
\caption{\label{fig_N2h3} Two nice sections with 4-crown stacks as retracts (hollow dots). The retraction of the 6-crown stack on the left is from \cite[Fig.\ 5.10]{Farley_1997}.}
\end{center}
\end{figure}

In a nice section $P$ not being a 2-antichain, the level-pairs $P(0,1)$ and $P(h_P-1,h_P)$ must form 6-crowns. The only nice sections of height 0, 1, and 2 are thus the 2-antichains, the 6-crowns, and the 6-crown stacks of height 2. All three are minimal automorphic. Two nice sections with height 3 are shown in Figure \ref{fig_N2h3}. Both contain a 4-crown stack of height 2 as retract.

\subsection{Crowned sections and split levels} \label{subsec_FourCrownStacks}

Our goal is a recursive approach for the determination of the nice sections of width three which have a 4-crown stack as retract. However, our results will concern a larger class of sections which contains the nice sections of width three as a proper subclass:

\begin{definitionaC} \label{def_crownedSection}
We call a section $P$ of width three a {\em crowned section} if $P(0,1)$ and $P(h_P-1,h_P)$ are 6-crowns. The symbol $\fC$ denotes the class of crowned sections.
\end{definitionaC} 

For a retraction $r : P \rarr R$ of $P \in  \fC$ onto a 4-crown stack, we are interested in proper lower or upper segments of $R$ which are the image of a proper lower or upper segment of $P$:

\begin{definitionaC} \label{def_splitLevel}
Let $P \in \fC$, let $r : P \rarr R$ a retraction onto a 4-crown stack, and let $R(\ell) = \setx{a,b}$ with $\lambda(a) \leq \lambda(b)$. We call $R(\ell)$ a {\em down-split level} if
$$
\ell < h_R \quad \mytext{and} \quad R(0 \rarr \ell) = r[ P(0 \rarr \lambda(b)) ],
$$
and we call $R(\ell)$ an {\em up-split level} if
$$
0 < \ell \quad \mytext{and} \quad R(\ell \rarr h_R) = r[ P(\lambda(a) \rarr h_R) ].
$$
\end{definitionaC}

In the retract of the poset on the left in Figure \ref{fig_N2h3}, the level set $R(0)$ is a down-split level, the level set $R(2)$ is an up-split level, and $R(1)$ is neither down- nor up-splitting. For the retract of the poset on the right of the figure, the down-split levels are $R(0)$ and $R(1)$ and the up-split levels are $R(1)$ and $R(2)$.

We agree on that for the rest of this section the poset $P$ is a crowned section and $r : P \rarr R$ is a retraction onto a 4-crown stack. This implies $h_R \geq 2$. Suppose to the contrary that $R$ is a 4-crown. For the corresponding spanning retraction $r' : P \rarr R'$, the poset $Q := P(0,1) \cup P(h_P)$ has the retraction $r' \vert_Q : Q \rarr R'$ onto the 4-crown $R'$. But this is impossible because $Q$ is a 6-crown or a 6-crown stack of height 2, or it contains a point $x \in Q(1)$ with $R'(0) < x < R'(1)$.

In our investigation of split levels, the level sets $R(0)$, $R(1)$, $R(h_R-1)$, and $R(h_R)$ will frequently play a special role. We notate the results for $R(0)$ and $R(1)$ only; the results for $R(h_R)$ and $R(h_R-1)$ are dual.

We start with a simple observation: Because the down-sets $r^{-1}(a)$, $a \in R(0)$, are pairwise disjoint,
\begin{equation} \label{R0CapP0_notEmpty}
R(0) \cap P(0) \not= \emptyset.
\end{equation}

\begin{lemmaaC} \label{lemma_Pkk}
Assume that the level set $R(\ell)$ is contained in a single level set of $P$. For $1 \leq \ell \leq h_R - 1$, the level set $R(\ell)$ is a split level, and for $\ell = 0$, one of the level sets $R(0)$ and $R(1)$ is a split level.
\end{lemmaaC}
\BP Let $R(\ell)$ be contained in $P(k)$ and let $z \in P(k) \setminus R(\ell)$.
\begin{itemize}
\item $\ell \in [1, h_R-1]$: Then $k \in [1, h_P-1]$. Depending on $\lambda(r(z)) \leq k$ or $\lambda(r(z)) \geq k$, the set $R(\ell)$ is a down-split level or an up-split level.
\item $\ell = 0$: Then $k = 0$ due to \eqref{R0CapP0_notEmpty}. If $r(z) \in R(0)$, then $R(0)$ is a down-split level. And in the case of $r(z) \notin R(0)$, we have $R(0) < x$ or $z < x$ for all $x \in P(1 \rarr h_P)$ because $P$ is crowned, and $R(1)$ is an up-split level.
\end{itemize}

\EP

Lemma \ref{lemma_R0P0} and Lemma \ref{lemma_Pkk} together yield that a crowned section has a 4-crown stack as retract iff it has a 4-crown stack as retract which contains a split level. On the one hand, this inconspicuous result brings our recursive approach into motion in Section \ref{subsec_recAppr}. On the other hand, this approach becomes more effective and the concept of split levels more applicable if we can show that in most cases we can expect additional split levels which are not induced by the action of $r$ on $P(0)$ and $P(h_P)$. We start with

\begin{lemmaaC} \label{lemma_Rell_maximalGestreut}
Let $R(\ell) = \setx{a,b}$ with $\lambda(a) < \lambda(b)$.
\begin{enumerate}
\item $\ell = 0$: The level set $R(1)$ is an up-split level if 
\begin{align*}
a & < P(  \min \lambda[ R(1) ]).
\end{align*}
\item $0 < \ell < h_R$: If 
\begin{align*}
a & < P(  \min \lambda[ R(\ell+1) ]) \\
\mytext{and} \quad b & > P(  \max \lambda[ R(\ell-1) ]),
\end{align*}
then $R(\ell-1)$ is a down-split level or $R(\ell+1)$ is an up-split level.
\end{enumerate}
\end{lemmaaC}
\BP 1) We have $\lambda(a) = 0$ according to \eqref{R0CapP0_notEmpty}. Due to $P(0) \setminx{a} < b$, we have $r(x) = b$ for all $x \in P(0) \setminx{a}$, thus $R(0) < r(y)$ for all $y \in P(\min \lambda[ R(1) ])$, and $R(1)$ is an up-split level.

2) Let $k := \lambda(a)$, $k^- := \max \lambda[ R(\ell-1)]$, $j :=  \lambda(b)$, and $j^+ := \min \lambda[ R(\ell+1)]$. Let $x \in P(k) \setminx{a}$. If $r(x) \in \setx{b} \cup R(\ell +1 \rarr h_R)$, then $P(k, j^+) \simeq 33$ yields $R(\ell) < r(z)$ for all $z \in P(j^+)$, and $R(\ell+1)$ is an up-split level.

Now assume $r(x) \in R(0 \rarr \ell - 1) \cup \setx{a}$ for all $x \in P(k) \setminx{a}$. Then $r(z) \in R(0 \rarr \ell - 1) \cup \setx{a}$ for all $z \in P(k^-)$, and $P(k^-, j) \simeq 33$ additionally yields $r(z) < b$, hence $r(z) \in R(0 \rarr \ell-1)$, and $R(\ell-1)$ is a down-split level.

\EP

The retracts in Figure \ref{fig_N2h3} illustrate the second part. In both, the level set $R(1)$ fulfills the assumptions, and $R(0)$ or $R(2)$ has to be a split level. In fact, both are.

\begin{figure}
\begin{center}
\includegraphics[trim = 75 715 425 75, clip]{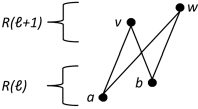}
\caption{\label{fig_abvw} The labeling of the points in $R(\ell,\ell+1)$ used in Lemma \ref{lemma_6Crown_33}, Corollary \ref{coro_6Cr_notSplitting}, and Proposition \ref{prop_doppel6Crown}. We assume $\lambda(a) \leq \lambda(b)$ and  $\lambda(v) \leq \lambda(w)$.}
\end{center}
\end{figure}

For the rest of this section, we deal with a 4-crown $R(\ell,\ell+1)$. As indicated in Figure \ref{fig_abvw}, we denote the points of $R(\ell)$ with $a$ and $b$ and assume $\lambda(a) \leq \lambda(b)$, and we denote the points of $R(\ell+1)$ by $v$ and $w$ and assume $\lambda(v) \leq \lambda(w)$.

\begin{lemmaaC} \label{lemma_6Crown_33}
Let $0 \leq \ell \leq h_R - 1$.

1) $P(\lambda(a),\lambda(w)) \simeq 33$.

2) If $P(\lambda(b),\lambda(v)) \simeq 33$, then $R(\ell)$ is a down-split level and $R(\ell+1)$ is an up-split level..

3) If $P(\lambda(b),\lambda(v))$ is a 6-crown and $P(\lambda(b),\lambda(w)) \simeq 33$, then $R(\ell)$ is a down-split level or $R(\ell+1)$ an up-split level.
\end{lemmaaC}
\BP 1) The points $a, \celaga{a}{b}$ and $w, \celaga{w}{v}$ form a 4-crown in $P(\lambda(a),\lambda(w))$.

2) The assumption yields $P(\lambda(b)) < R(\ell+1)$ and $R(\ell) < P(\lambda(v))$.

3) The point $v$ has two lower covers $b$ and $z$ in $P(\lambda(b)), \lambda(v))$, and we have $z < R(\ell+1)$, hence $r(z) \in R(0 \rarr \ell)$. For the point $y \in P(\lambda(b)) \setminx{b,z}$, we have $y < P(\lambda((v)) \setminx{v}$. If $r(y) \in R(0 \rarr \ell)$, then $R(\ell)$ is a down-split level, and in the case of $r(y) \in R(\ell+1 \rarr h_R)$, the level $R(\ell+1)$ is an up-split level.

\EP

\begin{figure}
\begin{center}
\includegraphics[trim = 75 685 325 75, clip]{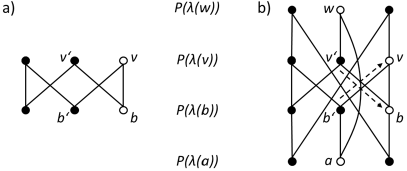}
\caption{\label{fig_6Cr_notSplitting} Illustration of Corollary \ref{coro_6Cr_notSplitting}.}
\end{center}
\end{figure}

\begin{corollaryaC} \label{coro_6Cr_notSplitting}
Let $1 \leq \ell \leq h_R - 2$ and let $P(\lambda(b),\lambda(v))$ form a 6-crown. Furthermore, as shown in Figure \ref{fig_6Cr_notSplitting}a, let $b'$ be the point in $P(\lambda(b)) \setminx{b}$ with $b' < v$, and let correspondingly $v'$ be the point in $P(\lambda(v)) \setminx {v}$ with $v' > b$.

If neither $R(\ell)$ nor $R(\ell+1)$ is splitting, then the poset induced by the points in $P(\lambda(a)) \cup P(\lambda(b)) \cup P(\lambda(v)) \cup P(\lambda(w))$ is isomorphic to the poset shown in Figure \ref{fig_6Cr_notSplitting}b. Furthermore, $r(v') = b$ and $r(b') = v$.
\end{corollaryaC}
\BP The relation $\lambda(a) < \lambda(b)$ is due to Lemma \ref{lemma_Pkk}. Because of the dual of Lemma \ref{lemma_6Crown_33}.3, the level pair $P(\lambda(a),\lambda(v))$ cannot be of type 33, and $P(\lambda(a),\lambda(v))$ must be a 6-crown. This implies $P(\lambda(a), \lambda(b)) \simeq 3C$, and $a < v$ yields $a < b'$. Due to Lemma \ref{lemma_6Crown_33}.1, the type of $P(\lambda(a), \lambda(w)))$ is 33. Finally, $r(v') \not= b$ would yield $r(v') \in R(\ell+1 \rarr h_R)$, and depending on the image of the point in $P(\lambda(v)) \setminx{v,v'}$, the level set $R(\ell)$ would be a down-split level or $R(\ell+1)$ an up-split level. The rest is dual.

\EP

\begin{propositionaC} \label{prop_doppel6Crown}
Let $R$ be a spanning retract and let $P(k,k+1)$ and $P(k+1,k+2)$ both be 6-crowns for a level index $k \in [1, h_P - 3]$. With
$$
\ell \; := \; \max \mysetdescr{ j \in [0,h_P] }{ \min \lambda[R(j)] \leq k },
$$
the level set $R(\ell)$ is a down-split level or the level set $R(\ell+1)$ is an up-split level.
\end{propositionaC}
\BP Because $R$ is spanning, the level index $\ell$ is well-defined with $\ell \leq h_R - 1$. With our convention about the points in $R(\ell,\ell+1) = \setx{a,b,v,w}$, the definition of $\ell$ implies $\lambda(a) \leq k$, $\lambda(b) \leq k+1$, $\lambda(v) \geq k+1$, and $P(\lambda(b),\lambda(v)) \not\simeq 3C$.

Assume $\ell = 0$. Then $\lambda(b) = 0$ and $\lambda(v) \geq k+1 \geq 2$ yield $P(\lambda(b),\lambda(v)) \simeq 33$, and Lemma \ref{lemma_6Crown_33}.2 delivers the result. The proof for $\ell+1 = h_R$ is dual.

Now let $1 \leq \ell \leq h_R - 2$. In the case of $P(\lambda(b),\lambda(v)) \simeq 33$, use Lemma \ref{lemma_6Crown_33}.2 again.

Suppose thus that $P(\lambda(b),\lambda(v))$ is a 6-crown. If the type of $P(\lambda(a),\lambda(v))$ or $P(\lambda(b),\lambda(w))$ is 33, then apply Lemma \ref{lemma_6Crown_33}.3 or its dual. Otherwise, we must have $\lambda(v) = k+1 = \lambda(w)$ and Lemma \ref{lemma_Pkk} delivers the result.

\EP

\subsection{A recursive approach} \label{subsec_recAppr}

The existence of a down-split level of $R$ ensures the existence of a proper lower segment $V$ of $P$ with a 2-antichain or a 4-crown stack as retract. Furthermore, with $D := \mysetdescr{ x \in P \setminus V }{ r(x) \in V }$, also the poset $P \setminus (V \cup D)$ has a 2-antichain or a 4-crown stack as retract. We give the constellation an own name:

\begin{definitionaC} \label{def_split}
Let $P \in \fC$. We call a quadrupel $(k, D, s, t)$ a {\em retractive down-split} of $P$ if
\begin{align*}
k & \in [0,h_P-1], \\
D &\subset P(k+1 \rarr h_P) \mytext{ is a down-set,} \\
s & : P(k+1 \rarr h_P) \setminus D \rarr S \mytext{ is a retraction,} \\
\mytext{and} \quad \quad \quad \quad 
t & : P(0 \rarr k) \rarr T \mytext{ is a retraction,} \\
\mytext{where } S \mytext{ and } T & \mytext{are 2-antichains or 4-crown stacks,}
\end{align*}
and we call a quadrupel $(k, U, s, t)$ a {\em retractive up-split} of $P$ if
\begin{align*}
k & \in [1,h_P], \\
U &\subset P(0 \rarr k-1) \mytext{ is an up-set,} \\
s & : P(0 \rarr k-1) \setminus U \rarr S \mytext{ is a retraction,} \\
\mytext{and} \quad \quad \quad \quad 
t & : P(k \rarr h_P) \rarr T \mytext{ is a retraction,} \\
\mytext{where } S \mytext{ and } T & \mytext{are 2-antichains or 4-crown stacks.}
\end{align*}
\end{definitionaC}

In order to avoid repetitions, we use the symbols $S$ and $T$ in what follows always as in this definition: for a retractive down-split $(k,D,s,t)$ of $P \in \fC$, we always have $S = s[P(k+1 \rarr h_P) \setminus D]$ and $T = t[ P(0 \rarr k)]$, and correspondingly for a retractive up-split. In a retractive down-split $(k,D,s,t)$, the poset $T$ is thus a retract of the lower segment of $P$ ending with level $k$ and $D$ is a down-set in $P(k+1 \rarr h_P)$, whereas in a retractive up-split $(k,U,s,t)$, the poset $T$ is a retract of the upper segment of $P$ starting with level $k$ and $U$ is an up-set in $P(0 \rarr k-1)$. The reader will observe that in a retractive down-split $(h_P-1,D,s,t)$ and a retractive up-split $(1,U,s,t)$ we must have $\# D \leq 1$ and $\# U \leq 1$, respectively, because $S$ contains at least two points.

We have seen in Section \ref{subsec_FourCrownStacks} that a crowned section has a 4-crown stack as retract iff it has such a retract with a split level. In terms of retractive splits this means

\begin{lemmaaC} \label{lemma_matching}
There exists a retraction of $P \in \fC_2$ onto a 4-crown stack iff 
\begin{itemize}
\item a retractive down-split $(k,D,s,t)$ of $P$ exists for which $T(h_T) < S(0)$ and for which the mapping
\begin{align*}
p & \mapsto
\begin{cases}
s(p), & \mytext{if } p \in P(k+1 \rarr h_P) \setminus D, \\
t(p), & \mytext{if } p \in P(0 \rarr k)
\end{cases}
\end{align*}
can be extended to a retraction $r : P \rarr T \oplus S$,
\item or a retractive up-split $(k,U,s,t)$ of $P$ exists for which $S(h_S) < T(0)$ and for which the mapping
\begin{align*}
p & \mapsto
\begin{cases}
s(p), & \mytext{if } p \in P(0 \rarr k-1) \setminus U, \\
t(p), & \mytext{if } p \in P(k \rarr h_P)
\end{cases}
\end{align*}
can be extended to a retraction $r : P \rarr S \oplus T$.
\end{itemize}
We call a retraction $r$ and a retractive split {\em matching} if they are coupled in this way.
\end{lemmaaC}

This result gives raise to a recursive approach to determine the crowned sections which have a 4-crown stack as retract. Let $\fS$ be a subclass of $\fC$ closed under duality, e.g., the subclass of nice sections of width three. For a retractive split of a poset $P \in \fS$, we need a {\em $t$-base} $V$ and an {\em $s$-base} $W$ of a poset $P \in \fS$ of width 3, i.e.,
\begin{itemize}
\item $V$ is a lower segment or an upper segment of $P$ which has a retraction $t$ onto a 2-antichain or a 4-crown stack $T$.
\item $W$ is an upper segment of $P$ which has a retraction of $W \setminus D$ onto a 2-antichain or a 4-crown stack $S$ for a down-set $D \subset W$, or $W$ is a lower segment of $P$ which has a retraction of $W \setminus U$ onto a 2-antichain or a 4-crown stack $S$ for an up-set $U \subset W$.
\item $V$ and $W$ are coupled by $V \cap W = \emptyset$ and $P = V \cup W$.
\end{itemize}

Assume that we have listed all lower segments up to height $n$ of sections contained in $\fS$ which have a 2-antichain or a 4-crown stack as retract. Given a poset $P \in \fS$ of height $n+1$, we can quickly identify all candidates for $t$-bases in $P$ by checking which of the lower segments in our list are lower segments $P(0 \rarr k)$ of $P$ and which of their duals are upper segments $P(k \rarr n+1)$ of $P$. Of course, $P(0)$ and $P(n+1)$ are always candidates. The poset $P$ has a 4-crown stack as retract iff one of the corresponding segments $P(k+1 \rarr n+1)$ or $P(0 \rarr k-1)$ provides an $s$-base which can be combined with the $t$-base as described in Lemma \ref{lemma_matching}.

In order to test the approach we need a sub-class of $\fC$ which combines plenty of interesting posets with good structural access. We choose the class of the nice sections with {\em horizon two}. We introduce them in the next section, and in Section \ref{sec_Application}, we identify all of them with height up to six with a 4-crown stack as retract.

\section{Nice sections with horizon two} \label{sec_horizonTwo}

\subsection{Definition and isomorphism types} \label{Nzwei_def}

For a crowned section $P \in \fC$ with height $h_P \geq 2$, we have $P(0,h_P) \simeq 33$. There exists thus a smallest integer $\eta$ with 
$$
P(k,k+\eta) \simeq 33 \quad \mytext{for all } k \in [0,h_P-\eta].
$$
We call this integer the {\em horizon} of $P$. It is the distance in which all details become blurred. The horizon of the sections in Figure \ref{fig_N2h3} is two. We have

\begin{propositionaC} \label{prop_C2istN2}
A crowned section with horizon 2 is nice.
\end{propositionaC}
\BP Let $P$ be a crowned section with horizon 2. We have to show that $P$ does not contain an irreducible point. Let $x \in P(k)$, $k \in [0, h_P-1]$. For the upper covers of $x$, only the case $P(k,k+1) \simeq 3C$ is of interest which yields $k \leq h_P - 2$. One upper cover of $x$ is contained in $P(k+1)$, and depending on the type of $P(k+1,k+2)$, the layer $P(k+2) > x$ contains one or even two additional upper covers of $x$. The argument for lower covers is dual.

\EP



Let $\fN_2$ denote the class of nice sections with horizon 2 and height at least two. The 6-crown stacks of height greater than one belong to $\fN_2$, and Farley \cite[Prop.\ 4.1]{Farley_1997} has shown that all 6-crown stacks with equal height are isomorphic.  The following theorem generalizes this result:

\begin{theoremaC} \label{theo_isomorphism}
The isomorphism type of a poset $P \in \fN_2$ is uniquely determined by the types of its level-pairs $P(k, k+1)$, $k \in [1, h_P-2]$. In particular, there exist $2^{n-2}$ isomorphism types of posets of height $n \geq 2$ in $\fN_2$. Furthermore, every permutation of $P(0)$ can be uniquely extended to an automorphism of $P$.
\end{theoremaC}
\BP Let $P = (X, \leq_P)$ and $Q = (Y,\leq_Q)$ be elements of $\fN_2$ and let $\alpha : P(0) \rarr Q(0)$ be a bijection. We define on $X$ and $Y$ the relations
\begin{align*}
\prec_P & := \mysetdescr{ (x,y) \in {<_P} }{ \lambda(y) - \lambda(x) = 1 } \\ 
\mytext{and} \quad \prec_Q & := \mysetdescr{ (x,y) \in {<_Q} }{ \upsilon(y) - \upsilon(x) = 1 },
\end{align*}
where $\upsilon$ is the level-function of $Q$. It is easily seen that we can extend $\alpha$ to an isomorphism $\beta : (X, \prec_P) \rarr (Y, \prec_Q)$ iff  $h_P = h_Q$ and
\begin{equation*}
\forall \; k \in [1,h_P-2] : \quad 
P(k,k+1) \mytext{ 6-crown} \; \; \Leftrightarrow \; \; Q(k,k+1) \mytext{ 6-crown,}
\end{equation*}
and that in this case $\beta$ is uniquely determined by $\alpha$. Every isomorphism from $(X, \prec_P)$ to $(Y, \prec_Q)$ sends $P(k)$ onto $Q(k)$ for all $k \in [0,h_P]$, and because $P$ and $Q$ have horizon 2, the isomorphisms between $(X, \prec_P)$ and $(Y, \prec_Q)$ are exactly the isomorphisms between $P$ and $Q$.

Now let $n \geq 2$ be an integer and let $K \subseteq \myNkz{n-1}$ be an index set with $0, n-1 \in K$. We construct a poset $P \in \fN_2$ with $h_P = n$ and $P(k, k+1)$ being a 6-crown iff $k \in K$. We start with the union $Q$ of three pairwise disjoint chains of length $n$, extend $Q(k,k+1)$ to an arbitrary 6-crown for every $k \in K$, and complete the partial order relation by adding $Q(k) \times Q(\ell)$ for all $k, \ell \in [0,n]$ with $k \leq \ell - 2$. It is easily seen that we can label the points in the resulting poset $P$ in such a way that all requirements on a section in Definition \ref{def_Section} are fulfilled. The poset $P$ is thus a crowned section with horizon 2, and because of Proposition \ref{prop_C2istN2}, it is nice.

\EP

Due to this result. we can always draw the Hasse-diagram of a poset $P \in \fN_2$ in a standardized way: all level-pairs $P(k,k+1)$ being 6-crowns are drawn with a ``void'' in the center as in Figure \ref{fig_N2h3}.

We denote the classes of segments, lower segments, and upper segments of nice sections with horizon 2 by $\fNS_2, \fNL_2$, and $\fNU_2$, respectively. Theorem \ref{theo_isomorphism} extends to these classes: Two segments of equal height $n$ are isomorphic iff both are 3-antichains (the case $n = 0$) or all of their consecutive level-pairs are isomorphic (the case $n \geq 1$), and also the extension of a level permutation to an automorphism is always possible. We will also use the following fact: if $t : V \rarr T$ is a retraction of a segment $V \in \fNS_2$ and $\alpha$ an automorphism of $V$, then $\alpha \circ t \circ \alpha^{-1}$ is a retraction of $V$ sharing all properties of $t$.




\subsection{Properties of retracts} \label{Nzwei_simpleResults}

We start our investigation with three simple observations:

\begin{corollaryaC} \label{coro_simpleResults}
Let $P \in \fNS_2$ and let $r : P \rarr R$ be a retraction.

1) For all $\ell \in [0,h_R]$ there exists an index $k \in [0,h_P-1]$ with $R(\ell) \subset P(k,k+1)$. In particular, $R(0) \subset P(0,1)$ if $\# R(0) \geq 2$.

2) Let $\ell \in [0,h_R]$ with $\# R(\ell) =2$ and $0 < \rho := \max \lambda[R(\ell)] < h_P$. We have $r[P(\rho+1 \rarr h_P)] = R(\ell+1 \rarr h_R)$ if
\begin{equation} \label{schubVonUnten}
r^{-1}(a) \cap P(\rho - 1) \not= \emptyset \quad \mytext{for all } a \in R(\ell).
\end{equation}
\end{corollaryaC}
\BP The first proposition follows from the definition of the horizon and \eqref{R0CapP0_notEmpty}. For the second one, observe that the assumption enforces $R(\ell) < r[P(\rho+1)]$.

\EP

\begin{figure}
\begin{center}
\includegraphics[trim = 70 725 445 75, clip]{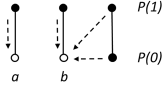}
\caption{\label{fig_P013C} The construction in the proof of Corollary \ref{coro_P01_3C}.}
\end{center}
\end{figure}

\begin{corollaryaC} \label{coro_P01_3C}
Let $P \in \fNS_2$ with $h_P \geq 2$ and $P(0,1) \simeq 3C$. P has a retract $R$ with $R(0,1)$ being a 4-crown iff $P(2 \rarr h_P)$ has a retract $W$ with $W(0)$ being a 2-antichain. 
\end{corollaryaC}
\BP ``$\Rightarrow$'': Let $R(0) = \setx{a,b}$ with $\lambda(a) \leq \lambda(b)$. Due to the first part of Corollary \ref{coro_simpleResults}, we have $R(0) \subset P(0,1)$. In the case of $\lambda(b) = 1$, the equation $r[P(2 \rarr h_P)] = R(1 \rarr h_P)$ follows with  \eqref{schubVonUnten}. In the case of $\lambda(b) = 0$, the relation $R(0) < P(2 \rarr h_P)$ enforces $r[P(2 \rarr h_P)] \subseteq R(1 \rarr h_P)$, and $P(0,1) \simeq 3C$ yields $R(1) \subset P(2 \rarr h_P)$.

``$\Leftarrow$'': If $w : P(2 \rarr h_P) \rarr W$ is a retraction with $W(0)$ being a 2-antichain, select two points $a, b \in P(0)$ and define a retraction $r : P \rarr W \cup \setx{a,b}$ by (cf.\ Figure \ref{fig_P013C})
\begin{align*}
r(x) & := 
\begin{cases}
w(x), & \mytext{if } x \in P(2 \rarr h_P), \\
a, & \mytext{if } x \in \setx{a, c_{1,\gamma(a)}}, \\
b, & \mytext{if } x \in P(0,1) \setminus \setx{a, c_{1,\gamma(a)}}.
\end{cases}
\end{align*}

\EP

\subsection{Retractive up- and down-splits} \label{Nzwei_retrSplits}

Nice sections with horizon two are good objects for our recursive approach because we can expect that a 4-crown stack $R$ as retract provides ``many'' split levels. Let $\ell \in [1,h_R-1]$ be a level index of $R$. Due to the horizon of $P$, the level set $R(\ell)$ is either contained in a single level set of $P$ or in two consecutive ones. Now the first statement in Lemma \ref{lemma_Pkk} and the second part of Lemma \ref{lemma_Rell_maximalGestreut} tell us that among the three consecutive level sets $R(\ell-1), R(\ell)$, and $R(\ell+1)$ there must be at least one split level. On the other hand, for a given retractive split of $P$ it is easy to check if there exists a matching retraction of $P$ onto a 4-crown stack:

\begin{theoremaC} \label{theo_RST}
Let $P \in \fN_2$. For a retractive down-split $(k,D,s,t)$ of $P$, the existence of a matching retraction $r : P \rarr T \oplus S$ is equivalent to
\begin{align} \label{forall_d_in_D}
\begin{split}
T(h_T) & <_P S(0),\\
\forall \; d \in D \;  \exists \; v \in T(h_T) :\; p & \not<_P d \mytext{ for all } p \in t^{-1}(v).
\end{split}
\end{align}
which implies $D \subseteq P(k+1)$, and for a retractive up-split $(k,U,s,t)$ of $P$, there exists a matching retraction $r : P \rarr S \oplus T$ iff
\begin{align} \label{forall_u_in_U}
\begin{split}
S(h_S) & <_P T(0),\\
\forall \; u \in U \;  \exists \; v \in T(0) :\; u & \not<_P p \mytext{ for all } p \in t^{-1}(v),
\end{split}
\end{align}
which implies $U \subseteq P(k-1)$.
\end{theoremaC}
\BP Let $r : P \rarr R$ be a retraction onto a 4-crown stack matching a retractive up-split $(k,U,s,t)$. Due to $r[P(k \rarr h_P)] = T$, the second part of \eqref{forall_u_in_U} is necessary for $r[U] \subseteq T$. The addendum follows because of $P(k-2) <_P P(k \rarr h_P) \supset T(0)$.

Now assume that \eqref{forall_u_in_U} holds for $(k,U,s,t)$ and let $u \in U \subseteq P(k-1)$. There exists a point $\tau(u) \in T(0)$ with $u \not<_P p$ for all $p \in t^{-1}(\tau(u))$. We define $\rho(u)$ as the single point contained in $T(0) \setminus \setx{ \tau(u) }$. We define a mapping $r : P \rarr P$ by setting for $x \in P$
\begin{align*}
r(x) & :=
\begin{cases}
s(x), & \mytext{if } x \in P(0 \rarr k-1) \setminus U, \\
\rho(x) , & \mytext{if } x \in U, \\
t(x), & \mytext{if } x \in P(k \rarr h_P),
\end{cases}
\end{align*}
The mapping $r$ is clearly idempotent. In order to show that it is order-preserving, let $x <_P y$. Only three cases have to be discussed:
\begin{itemize}
\item $x \in P(0 \rarr k-1) \setminus U$ and $y \in U \cup P(k \rarr h_P)$ yields $r(x) \in S <_{S \oplus T} T \ni r(y)$.
\item $x \in U$, $y \in P(0 \rarr k-1)$ is not possible due to $U \subseteq P(k-1)$.
\item $x \in U$ and $y \in P(k \rarr h_P)$: $x \not< p$ for all $p \in t^{-1}(\tau(x))$ delivers $t(y) \in T \setminx{\tau(x)}$. Because $\rho(x)$ is the only minimal point of this poset, $r(x) = \rho(x) \leq t(y) = r(y)$.
\end{itemize}

\EP

\begin{figure}
\begin{center}
\includegraphics[trim = 75 720 335 75, clip]{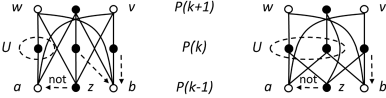}
\caption{\label{fig_CoroLuecke} Illustrations for the proof of Corollary \ref{coro_stackMitLuecke}. The point $z$ is not mapped to $a$ by the retraction $s$ of $P(0 \rarr k-1)$, and the points in $P(k) \setminus U$ are mapped to $b$ by the retraction $s'$ of $P(0 \rarr k) \setminus U$.}
\end{center}
\end{figure}

\begin{corollaryaC} \label{coro_stackMitLuecke}
Let $k \in [1,h_P-1]$ and assume that we have retractions $s : P(0 \rarr k-1) \rarr S$ and $t : P(k+1 \rarr h_P) \rarr T$ with $S$ and $T$ being 2-antichains or 4-crown stacks with
\begin{align*}
s^{-1}(a) & = \setx{a} \mytext{ for a point } a \in S(h_S), \\
\mytext{and} \quad t^{-1}(v) & = \setx{v} \mytext{ for a point } v \in T(0).
\end{align*}
$P$ has a retract isomorphic to $S \oplus T$ if $P(k-1 \rarr k+1)$ is not a 6-crown stack.
\end{corollaryaC}

\BP Let $S(h_S) = \setx{a,b}$ and $T(0) = \setx{v,w}$. The assumptions imply $a \in P(k-1)$ and $v \in P(k+1)$. Furthermore, $s[P(k-1)\setminx{a}] \subseteq S \setminx{a}$ and $t[P(k+1)\setminx{v}] \subseteq T \setminx{v}$.

Assume $P(k,k+1) \simeq 3C$ (the other case is dual) and let $U \subset P(k)$ be the set of upper covers of $a$ in $P(k)$. Applying Theorem \ref{theo_isomorphism} on $P(k+1 \rarr h_P)$, we see that without loss of generality we can assume that none of the points in $U$ is below $v$. Figure \ref{fig_CoroLuecke} illustrates the case $R(\ell) \subset P(k-1)$, $R(\ell+1) \subset P(k+1)$.

No point of $P(k) \setminus U$ is above $a$ and the point $b$ is the only maximal point of $S \setminx{a}$. Due to $s^{-1}(a) = \setx{a}$, the mapping $s' : P(0 \rarr k) \setminus U \rarr S$ defined by
\begin{align*}
s'(x) & :=
\begin{cases}
s(x), & \mytext{if } x \in P(0 \rarr k-1), \\
b, & \mytext{if } x \in P(k) \setminus U
\end{cases}
\end{align*}
is a retraction, and due to $t^{-1}(v) = \setx{v}$, the quadrupel $(k+1, U, s', t)$ is a retractive up-split fulfilling \eqref{forall_u_in_U}.

\EP

\section{Application} \label{sec_Application}

With the recursive approach described in Section \ref{subsec_recAppr}, we will now identify all posets $P \in \fNL_2$ with height up to six which have a 4-crown stack as retract. For some of them, e.g. the 6-crown stacks, we already know the result. We include them in our investigation because we want to test our approach with as many nice sections as possible.

We encode the lower segments not being 3-antichains by the sequence of their level-pair types (cf.\ Theorem \ref{theo_isomorphism}). The letter ``1'' indicates a 6-crown, ``0'' a type $3C$. The lower segment 111 is thus the 6-crown stack of height 3, whereas 110 indicates the lower segment starting with two 6-crowns followed by a single $3C$.

\begin{table}
\begin{center} 
{\scriptsize
\begin{tabular}{| l l c | l l c | l l c | l l c | }
\hline
1 & & n & & & & & & & & &  \\
\hline
11 & & n & 10 & & y & & & & & & \\
\hline
111 & & y & 101 & & y & 110 & & n & 100 & & n \\
\hline
1111 & 0,3 & n & 1101 & 0 & n & 1011 & 0,2,3 & n & 1001 & 0,2 & y \\
1110 & & n & 1100 & & n & 1010 & & y & 1000 & & y \\
\hline
11111 & 0,3 & n & 11101 & 0,3 & y & 11011 & 0 & n & 10111 & 0,2,3 &  y \\
11001 & 0 & y & 10101 & 0,2,3,4 & y & 10011 & 0,2,4 & y & 10001 & 0,2,4 & y \\
11110 & & y & 11100 & & y & 11010 & & n & 10110 & & y \\
11000 & & n & 10100 & & y & 10010 & & n & 10000 & & n \\
\hline
111111 & 0,3 & y & 111101 & 0,3,5 & y & 111011 & 0,3,5 & n & 110111 & 0 & n \\
101111 & 0,2,3,5 & y & 111001 & 0,3,5 & n & 110101 & 0 & n & 101101 & 0,2,3,5 & y \\
110011 & 0,5 & n & 101011 & 0,2,3,4,5 & n & 100111 & 0,2,4,5 & n & 110001 & 0 & n \\
101001 & 0,2,3,4,5 & y & 100101 & 0,2,4 & y & 100011 & 0,2,4,5 & n & 100001 & 0,2,4 & y \\
111110 & & n & 111100 & & n & 111010 & & n & 110110 & & n \\
101110 & & n & 111000 & & n & 110100 & & n & 101100 & & n \\
110010 & & n & 101010 & & y & 100110 & & y & 110000 & & n \\
101000 & & y & 100100 & & y & 100010 & & y & 100000 & & y \\
\hline
\end{tabular} }
\caption{\label{table_LS6} The lower segments in $\fNL_2$ with height from one to six. Explanation in text.}
\end{center}
\end{table}

\begin{figure}
\begin{center}
\includegraphics[trim = 75 695 195 75, clip]{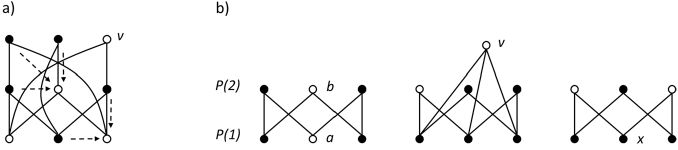}
\caption{\label{fig_Retr10Crit5} a) The lower segment $10$ and a retraction $r$ onto a 4-crown with $r^{-1}(v) = \setx{v}$ for a point $v$ from the top level of the 4-crown. b) Illustrations for the case discrimination in the proof of Criterion \ref{crit_split}.2. The hollow dots belong to $S(0)$.}
\end{center}
\end{figure}

The lower segments with height from one to six are listed in Table  \ref{table_LS6}. The letters ``y'' and ``n'' indicate the result of our investigation whether the segment has a 4-crown stack as retract or not. For the posets 1, 10, 11, 111, and 101, we refer to common knowledge. The lower segment $10$ even has a retraction $r$ to a 4-crown with $r^{-1}(v) = \setx{v}$ for a point $v$ from the top level of the 4-crown, as shown in Figure \ref{fig_Retr10Crit5}a.

For the posets with a final 0, we can apply the dual of Corollary \ref{coro_P01_3C} on previous results in Table \ref{table_LS6}. As an example, the posets $110$ and $100$ do not have a 4-crown stack as retract, because the 6-crown $P(0,1)$ does not have a suitable retract. And $1010$ and $1000$ have a 4-crown stack as retract because 10 is marked with ``y''.

The bulk of the work concerns thus the lower segments $P$ with a 6-crown as final level-pair. They all are elements of $\fN_2$. For them, the integers in Table \ref{table_LS6} indicate the levels $k \in [0, h_P-1]$ for which $P(0 \rarr k)$ can be a $t$-base according to the results for lower segments with less height. For $1111$ it is 0 and 3: 0, because the antichain $P(0)$ is always a candidate for a $t$-base, and 3, because the poset $1111(0 \rarr 3) = 111$ is marked with ``y'' in the previous row. Because 1 and 11 are marked with ``n'', the level indices 1 and 2 are missing in the list for $1111$.

The following criteria will be useful in our investigation:

\begin{criteriaaC} \label{crit_split}
Let $P \in \fN_2$ with $h_P \geq 3$. Assume that there exists a matching retraction $r$ for the retractive down-split $(k,D,s,t)$.
\begin{enumerate}
\item $k=0$ implies that $D \subset P(1)$ contains at most a single point.
\item $k=0$ and $P(1,2)$ being a 6-crown imply together that the poset $P(3 \rarr h_P)$ has a 2-antichain or a 4-crown stack as retract.
\item $k = h_P-2$ implies that $D \subseteq P(k-1)$ contains at least two points.
\item $k = h_P-2$ and $P(h_P-2,h_P-1)$ being a 6-crown imply together that the poset $P(0 \rarr h_P-3)$ has a 2-antichain or a 4-crown stack as retract.
\item $k=h_P-1$ implies that the poset $P(0 \rarr h_P-3)$ has a 2-antichain or a 4-crown stack as retract.
\end{enumerate}
\end{criteriaaC}
\BP 1. Trivial.

2. We discriminate three cases (cf.\ Figure \ref{fig_Retr10Crit5}b). If $S(0)$ contains a point $a \in P(1)$, then it has to contain a point $b \in P(2)$, too, because of the first condition in \eqref{forall_d_in_D}. At least one of the points in $P(1) \setminus \setx{a}$ has to be mapped to $b$ (first criterion), thus $s[P(3 \rarr h_P)] = S(1 \rarr h_S)$ due to \eqref{schubVonUnten}.

Now assume that $S(0)$ does not contain a point of $P(1)$. If $S(0)$ contains a point $v \in P(3 \rarr h_P)$, then at least two points in $P(1)$ are below $S(0)$ in contradiction to the first criterion. And in the case of $S(0) \subset P(2)$, exactly one point $x \in P(1)$ is below $S(0)$ and has to be sent to $P(0)$. From the two remaining points in $P(1)$ each is under a single different point of $S(0)$. Due to the first criterion, they have to be mapped onto $S(0)$, and $s[P(3 \rarr h_P)] = S(1 \rarr h_S)$ follows again with \eqref{schubVonUnten}.

3. The poset $P(h_P-1,h_P) \setminus D$ has to be disconnected.

4. According to the third criterion, the set $D \subset P(h_P-1)$ contains at least two points. The second condition in \eqref{forall_d_in_D} says that every point in $D$ is comparable with a single point of $T(h_T)$ only which implies $T(h_T) \subset P(h_P-2)$. A matching retraction maps $D$ onto $T(h_T)$, and the dual of \eqref{schubVonUnten} yields the result.

5. The first condition in \eqref{forall_d_in_D} yields $T(h_T) \not\subset P(h_P-1)$. Now the dual of \eqref{schubVonUnten} yields $t[P(0 \rarr h_P-3)] = T(0 \rarr h_T - 1)$.

\EP

\begin{figure}
\begin{center}
\includegraphics[trim = 70 670 240 70, clip]{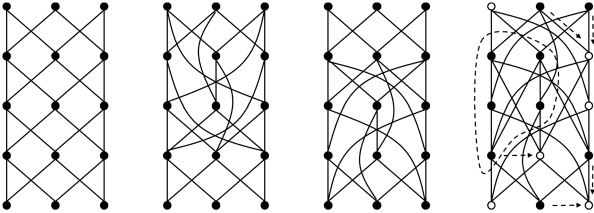}
\caption{\label{fig_PosetsHeight4} The posets $1111$, $1101$, $1011$, and $1001$. As already mentioned at the end of the introduction, the poset $1011$ is also treated in \cite[p.\ 99]{Schroeder_2016} and mentioned in \cite[p.\ 136]{Farley_1997} as a very complicated case. For the poset $1001$, a retraction onto a 4-crown tower is indicated which results from a retractive down-split $(2,D,s,t)$ based on the retraction $t$ of $P(0 \rarr 2) = 10$ shown in Figure \ref{fig_Retr10Crit5}a.}
\end{center}
\end{figure}

The rightmost poset in Figure \ref{fig_PosetsHeight4} confirms that we can neither drop the condition ``$P(1,2)$ is a 6-crown'' in Criterion \ref{crit_split}.2 nor the condition ``$P(h_P-2,h_P-1)$ is a 6-crown'' in Criterion \ref{crit_split}.4.

\begin{figure}
\begin{center}
\includegraphics[trim = 75 670 210 75, clip]{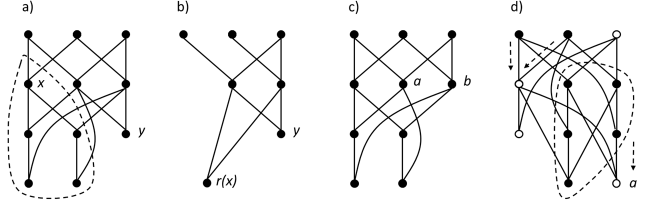}
\caption{\label{fig_011_001} The first three diagrams are illustrations for the proof of Lemma \ref{lemma_Q011}. Additionally the rightmost diagram shows for $P = 001$ a retraction of $P \setminus D$ to a 4-crown with a singleton $D \subset P(0)$.}
\end{center}
\end{figure}

\begin{lemmaaC} \label{lemma_Q011}
Let $P = 011$. If $r : P \setminus D \rarr R$ is a retraction onto a 4-crown stack or a 2-antichain for a down-set $D$, then $D$ contains at least two points of $P(0)$.
\end{lemmaaC}

\BP Assume that a retraction exists for $D = \setx{d}$ with $d \in P(0)$. Figure \ref{fig_011_001}a shows $P \setminx{d}$. Because $r$ cannot be extended to $P$, both upper covers $x$ and $y$ of $d$ have to be mapped to different points of $R(0)$. The encircled set $\darr x$ has thus to be mapped to a single point in $R(0)$. The resulting poset is shown in Figure \ref{fig_011_001}b. Clearly, it does not have a 4-crown stack or a 2-antichain as retract.

Now assume $D = \setx{d,y}$ with $d \in P(0)$ and $d \lessdot y$. The poset $P \setminus D$ is shown in Figure \ref{fig_011_001}c. If a retraction as described in the lemma exists, it cannot be extended to $P \setminx{d}$. The upper covers $a$ and $b$ of $y$ are thus mapped to different points of $R(0)$. But this is impossible because of $P(0) \setminx{d} < \setx{a,b}$.

\EP

Now we can start with the posets of height four. They are shown in Figure \ref{fig_PosetsHeight4}.

\begin{itemize}
\item P = 1111: Criterion \ref{crit_split}.5 prevents $P(4)$ from being a candidate for an $s$-base because the 6-crown $P(0,1)$ does not have a suitable retract. Furthermore, selecting $P(1 \rarr 4)$ as $s$-base does not work due to Criterion \ref{crit_split}.2, because $P(3,4) \simeq 1$ is marked with ``n'' in Table \ref{table_LS6}. $1111$ does thus have no retractive down-split fulfilling \eqref{forall_d_in_D}, and because the poset is self-dual, we conclude that it does not have a 4-crown stack as retract.
\item $P = 1101$: Again, $P(1 \rarr 4)$ fails as an $s$-base due to Criterion \ref{crit_split}.2. Therefore, $1101$ does not have a retractive down-split fulfilling \eqref{forall_d_in_D}.
\item $P = 1011$: $P(1 \rarr 4)$, $P(3,4)$, or $P(4)$ being a successful $s$-base is prevented by Criterion \ref{crit_split}.1 together with Lemma \ref{lemma_Q011}, Criterion \ref{crit_split}.4, and Criterion \ref{crit_split}.5, respectively. Also $1011$ does thus not have a retractive down-split fulfilling \eqref{forall_d_in_D}. Therefore, neither $1101$ nor $1011$ has a 4-crown stack as retract.

\item $P = 1001$: A retraction to a 4-crown stack resulting from a retractive down-split $(2,D,s,t)$ is indicated in Figure \ref{fig_PosetsHeight4} on the right.
\end{itemize}

We continue with the posets of height five:

\begin{itemize}
\item $P = 11111$: Due to the Criteria \ref{crit_split}.2 and \ref{crit_split}.4, the upper segments $P(1 \rarr 5)$ and $P(4,5)$ cannot work as $s$-bases. The self-dual poset has thus no 4-crown stack as retract.
\item $P = 11101$: see $P = 10111$.
\item $P = 11011$: Criterion \ref{crit_split}.2 prevents $P(1 \rarr 5)$ from being a successful $s$-base. Because $11011$ is self-dual, the poset cannot have a 4-crown stack as retract.
\item $P = 10111$: Figure \ref{fig_N2h3} shows that the 6-crown stack $Q : P(2 \rarr 5) = 111$ has a retraction $t : Q \rarr T$ with $T$ being a 4-crown stack and $t^{-1}(a) = \setx{a}$ for a point $a \in T(0)$. With $s : P(0) \rarr P(0)$ being a mapping such that $s[P(0)]$ is a 2-antichain, apply Corollary \ref{coro_stackMitLuecke} with $k = 1$.
\item $P = 11001$: see $10011$.
\item $P = 10101, 10001$: Let $s : P(0 \rarr 2) \rarr S$ and $t : P(3 \rarr 5) \rarr T$ be retractions onto 4-crowns as shown in Figure \ref{fig_Retr10Crit5}a and its dual. Without loss of generality (Theorem \ref{theo_isomorphism}) we assume that the single point in $P(2) \cap S$ is below the single point in $P(3) \cap T$. Now $(3, \emptyset, s, t)$ is a retractive down-split fulfilling \eqref{forall_d_in_D}.

\item $P = 10011$: Let $Q := 1001$. The retraction $t $ of $Q$ indicated in Figure \ref{fig_PosetsHeight4} on the right and a self-mapping $s$ of $P(5)$ onto a 2-antichain can be combined to a retractive down-split $(4, \emptyset, s, t)$ of $P$ fulfilling \eqref{forall_d_in_D}.
\end{itemize}

For some of the posets of height six, a quick positive decision is possible:

\begin{itemize}
\item $P = 111111, 111101, 101111, 101101$: Combine the retractions in Figure \ref{fig_N2h3}.
\item $P = 101001, 100001$: Let $P$ be one of these posets. With $Q := P(3 \rarr 6) = 001$, select a point $d \in Q(0)$ and let $s$ be the retraction of $Q \setminus \setx{d}$ onto a 4-crown $S$ shown in Figure \ref{fig_011_001}d with $a \in S \cap Q(0)$ as in the figure. Furthermore, let $t$ be the retraction of $P(0 \rarr 2) = 10$ onto a 4-crown $T$ as shown in Figure \ref{fig_Retr10Crit5}a with $v \in T \cap P(2)$ as in the figure. Applying Theorem \ref{theo_isomorphism} on $P(0 \rarr 2)$, we can assume $v \lessdot a$ and $v \not< d$ and get a retractive down-split $(2,\setx{d},s,t)$ fulfilling \eqref{forall_d_in_D}.
\item $P = 100101$: dual to $101001$.
\end{itemize}

For $P = 110011$, the Criteria \ref{crit_split}.2 and \ref{crit_split}.5 prevent $P(1 \rarr 6)$ and $P(6)$ from being successful $s$-bases; the self-dual poset has thus now 4-crown stack as retract. The remaining posets of height six are $111011$, $111001$, $101011$, $110001$, and their duals. They do not have a 4-crown stack as retract:

\begin{lemmaaC} \label{lemma_P36_gehtNicht}
None of the posets $111011$, $111001$, $101011$, and $110001$ has a 4-crown stack as retract.
\end{lemmaaC}
\BP Firstly, let $P \in \fN_2$ be any nice section of height six and $R$ a retract of $P$ which is a 4-crown stack. Every level set of $R$ has to be contained in a single level set or in two consecutive level sets of $P$, and we conclude that at least one of the posets $P(0 \rarr 3) \cap R$ and $P(3 \rarr 6) \cap R$ is a 4-crown stack.

Now let $P$ be any of the four posets and $r$ a retraction onto a 4-crown stack $R$. With $Q := P(3 \rarr 6) \cap R$, we show in a first step that $Q$ cannot be a 4-crown stack.

$P(0 \rarr 3) = 111, 110$: $P(3 \rarr 6)$ does not have a 4-crown stack as retract. $Q$ being a 4-crown stack requires thus a point $x \in P(3)$ with $v := r(x) \in P(0 \rarr 2) = 11$. Let $\setx{v,w}$ be the level set of $R$ containing $v$. $x \in P(3)$ yields $r[P(0,1)] \leq v$, hence $w \notin P(0,1)$. But then, $r(z) < \setx{v,w}$ holds for at least five of the points $z \in P(0,1)$, and $P(0,1)$ cannot be mapped onto a 2-antichain.

$P = 101011$: $Q$ being a 4-crown stack requires due to Lemma \ref{lemma_Q011} that two points $x, y \in P(3)$ are mapped into $P(0 \rarr 2) = 10$. The points $r(x)$ and $r(y)$ belong thus to $R(0,1)$. Three cases are possible, each leads to a contradiction:
\begin{itemize}
\item $\setx{r(x),r(y)}$ is a level set of $R$: Let $z \in P(2)$ be the common lower cover of $x$ and $y$. Due to $P(0) < z$, the point $r(z)$ is the only minimal point of $R$.
\item $r(x) < r(y)$: Then $r(x) \in R(0)$, and $P(0) < x$ again yields that $r(x)$ is the only minimal point of $R$.
\item $v := r(x) = r(y)$: Let $\setx{v,w}$ be the level set of $R$ containing $v$. Due to $r[P(0 \rarr 2)] \leq v$, the point $w$ must belong to $P(3 \rarr h_P)$, and due to $\lambda(v) \leq 2$, we must have $w \in P(3)$ and $v \in P(2)$. But then $w \in Q \not\ni v$.
\end{itemize}

$Q$ is thus not a 4-crown stack. This requires a level set $R(\ell) = \setx{a,b}$ with $a \in P(2)$ and $b \in P(3)$. Clearly, $\ell = 1$. Due to $P(1) < b$, no point of $P(1)$ can be mapped to $a$. And if a point of $P(1)$ is mapped to $b$, \eqref{schubVonUnten} yields $R(2 \rarr h_R)$ being a retract of $P(4 \rarr 6) = 11$ which is impossible. But if no point of $P(1)$ can be mapped to $R(1)$, then $r[P(0,1)] = R(0)$ which is impossible, too.

\EP








\end{document}